\documentclass[reqno,12pt]{amsart}
\newtheorem{theo}{Theorem}
\newtheorem{rem}{Remark}

\newtheorem{prop}{Proposition}

\newcommand\eps\varepsilon
\newcommand\ph\varphi
\newcommand\kap\varkappa

\usepackage{enumerate}

\begin{document}

\title[ Coincidence points of Mappings in Banach Spaces]
{ Coincidence points of Mappings in Banach Spaces }

\author[Oleg Zubelevich]{Oleg Zubelevich\\ \\\tt
 Dept. of Theoretical mechanics,  \\
Mechanics and Mathematics Faculty,\\
M. V. Lomonosov moscow State University\\
Russia, 119899, Moscow,  MGU \\
 }
\date{}
\thanks{Partially supported by grants
 RFBR  18-01-00887}

\begin{abstract}In this article we prove an existence theorem for coincidence points of  mappings in Banach spaces. This theorem generalizes the Kantorovich fixed point theorem.
\end{abstract}

\maketitle

\numberwithin{equation}{section}
\newtheorem{theorem}{Theorem}[section]
\newtheorem{lemma}[theorem]{Lemma}
\newtheorem{definition}{Definition}[section]

\section{Introduction}

In this short note we develop the ideas of Leonid Kantorovich \cite{KA} to the problem of coincidence points of mappings in Banach spaces. 

Problem of coincidence points is  actively discussed for the last decade. We do not  try to present any survey of the large amount of articles   published in this field.

We just mention only one very general and simple result obtained by A.V. Arutyunov that is most close to the topic of our  article. 

Let $(Q,d_Q),(P,d_P)$ stand for the metric spaces; the space $Q$ is complete. And let  $u,v:Q\to P$ be mappings.

By $v:Q\to P$  denote a Lipschitz continuous mapping: 
$$d_P(v(x_1),v(x_2))\le \beta\, d_Q(x_1,x_2),\quad x_1,x_2\in Q.$$
A mapping $u:Q\to P$ we assume to be continuous and such that the inequality
$$d_P(u(x'),y)\le \alpha\,r$$ implies that there exists  an element $x\in Q,\quad d_Q(x,x')\le r$ which satisfies the following equation $u(x)=y$.
This assumption holds for any $r>0$; the positive number $\alpha$ is fixed. 

Such a mapping $u$ is said to be the $\alpha-$covering.

\begin{theo}[\cite{arut}]\label{sdfddd} If $\beta<\alpha$ then  there exists a point $x_*\in Q$ such that $v(x_*)=u(x_*).$\end{theo}
The proof of this theorem is based upon a successive approximation procedure that is similar to one from the contraction mapping principle. 

Describe briefly this process. Take any point $x_0\in Q$. Then there exists a point $x_1\in Q$ such that $u(x_1)=v(x_0)$ and so on. Eventually, we obtain $u(x_{i+1})=v(x_i)$. The sequence  $\{x_i\}$ can be chosen such that $$d_Q(x_{i+1},x_i)\le \frac{\beta}{\alpha}d_Q(x_{i},x_{i-1}).$$ So it is convergent with exponential rate.

In this article we present a theorem that allows to obtain existence results even for the case when the successive approximation procedure converges arbitrary slow.

\subsection{Example}

Let $X,Y$ be Banach spaces over the field $\mathbb{R}$ or $\mathbb{C}$. 
By 
$B_X(\hat x,R)$ we denote the  open ball of this space:
$$B_X(\hat x,R)=\{x\in X\mid \|x-\hat x\|_X<R\},$$ and $ \overline B_X(\hat x,R)=\{x\in X\mid \|x-\hat x\|_X\le R\}.$

Consider a continuous  bilinear function $A:X\times X\to Y$ such that
$$\|A(x_1,x_2)\|_Y\le a\|x_1\|_X\|x_2\|_X,\quad A(x_1,x_2)=A(x_2,x_1).$$
Let $B:X\to Y$ stand for a bounded linear operator onto: 
$$\overline B_Y(0,b)\subseteq B(\overline B_X(0,1))$$ with some positive constant $b$. Let $C$ stand for a fixed element of $Y,\quad \|C\|_Y=c$. 
\begin{prop}\label{dgttt}If $D=b^2-4ac\ge 0$ then the equation
\begin{equation}\label{dgyy}A(x,x)+Bx+C=0\end{equation} has a solution.\end{prop}
If $D>0$ then this proposition  follows from theorem \ref{sdfddd}; the corresponding successive approximations converge exponentially. If $D=0$ then the successive procedure converges slower than exponentially and the result does not follow from theorem \ref{sdfddd}.

Proposition \ref{dgttt} follows from theorem \ref{sd4te}, we discuss it in the next section.

\section{Main Theorem}  

Suppose for some fixed $x_0\in X$ and $r>0$ we are given with a continuous mapping 
$$\Psi: B_X(x_0,2r)\to Y.$$

In the sequel we use a notation $I=[\tau_0,\tau_0+r)$; it is allowed to be $r=\infty$. Let $\psi\in C(I)$ be a non-decreasing  function.

\begin{definition}\label{csf} We shall say that $\Psi$ is a $\psi$-covering iff for any $\tau_0\le\tau'<\tau''< \tau_0+r$ and for any $x'\in B_X(x_0,r)$ one has
$$ \overline B_Y(\Psi(x'),\psi(\tau'')-\psi(\tau'))\subseteq \Psi\Big(  \overline B_X(x',\tau''-\tau')\Big).$$
\end{definition}

\begin{rem}
If we take $\psi(\tau)=\alpha\tau$ then $\psi-$covering turns into $\alpha-$covering.\end{rem}

Introduce also a function
$\Phi:B_X(x_0,r)\to Y$ which is assumed to be 
 Frech\'et  differentiable at each point and the function $x\mapsto \Phi'(x)$ is continuous with respect to the operator norm.

\begin{theo}\label{sd4te}
Suppose that the following hypotheses hold
\begin{enumerate}
    \item[H1.]\label{h1} The mapping $\Psi$ is a $\psi-$covering. 
\item[H2.]\label{h2} There exists an  increasing function $\ph\in C^1(I)$ such that  $$\|\Phi(x_0)-\Psi(x_0)\|_Y\le \ph(\tau_0)-\psi(\tau_0)$$ and
for each $\tau\in I$ the inequality $\|x-x_0\|_X\le\tau-\tau_0$ implies $$\|\Phi '(x)\|\le \ph'(\tau).$$\end{enumerate}

Assume also that an equation $\psi(\tau)=\ph(\tau)$ has at least one solution and let $\tau_*\in I$ be its smallest solution. Then there exists a point $$x_*\in B_X(x_0,r),\quad \|x_*-x_0\|_X\le \tau_*-\tau_0$$ such that
\begin{equation}\Phi(x_*)=\Psi(x_*).\label{sv}\end{equation}
\end{theo}
\begin{rem} If $X=Y,\quad \Psi(x)=x,\quad \psi(\tau)=\tau$ then hypothesis H1 holds automatically and theorem \ref{sd4te} turns into Kantorovich's result.\end{rem}

Rewrite equation (\ref{dgyy}) as follows
$$\Phi(x)=A(x,x)+C,\quad \Psi(x)=-Bx,\quad \Phi(x)=\Psi(x).$$
It is not hard to show that $\Psi$ is a $\psi-$covering with $$\psi(\tau)=b\tau,\quad r=\infty,\quad \tau_0=0,\quad x_0=0$$ and $\ph(\tau)=a\tau^2+c.$

 \section{Proof of Theorem \ref{sd4te}}
Show that there is an increasing sequence $\{\tau_j\}_{j\in\mathbb{N}}\subset I$ such that $\tau_j\to \tau_*$ and $\psi(\tau_{j+1})=\ph(\tau_j)$.  

Indeed, introduce a function $\alpha_0(\tau)=\psi(\tau)-\ph(\tau_0)$. Hypothesis H2 implies $\alpha_0(\tau_0)<0$ but since $\ph$ is an increased function we see $$\alpha_0(\tau_*)=\psi(\tau_*)-\ph(\tau_0)>0.$$ Thus there is $\tau_1\in (\tau_0,\tau_*),\quad \alpha_0(\tau_1)=0.$

Take a function $\alpha_1(\tau)=\psi(\tau)-\ph(\tau_1)$. It is easy to see
$$\alpha_1(\tau_*)>0,\quad \alpha_1(\tau_1)<0$$ thus there exists $\tau_2\in(\tau_1,\tau_*)$ such that $\alpha_1(\tau_2)=0$ and so on.

\begin{lemma}\label{sdfsrrrrr}
There exists a sequence $\{x_j\}\subset X,\quad j=0,1,2,\ldots$ such that
\begin{align}
\|x_j-x_0\|_X&\le \tau_j-\tau_0,\label{l1}\\
\|x_{j+1}-x_{j}\|_X&\le \tau_{j+1}-\tau_{j},\label{l2}\\
\Psi(x_{j+1})&=\Phi(x_{j}),\label{l3}
\end{align}\end{lemma}
\begin{rem}From (\ref{l1}) and (\ref{l2})   it follows that $x_j,x_{j+1}\in B_X(x_0,r)$.\end{rem}

{\it Proof of Lemma \ref{sdfsrrrrr}} We build this sequence by induction. Let us obtain formulas (\ref{l1})-(\ref{l3}) for $j=0$.

To employ definition \ref{csf} we put $$\tau'=\tau_0,\quad x'=x_0,\quad \tau''=\tau_1,\quad y=\Phi(x_0).$$ 
Then in accordance to hypothesis H2 $$\|y-\Psi(x')\|_Y\le \psi(\tau'')-\psi(\tau').$$ and by hypothesis H1 there exists $x_1,\quad \|x_1-x_0\|_X\le \tau_1-\tau_0$ such that $$\Psi(x_1)=\Phi(x_0).$$

Assume that we have already constructed $x_1,\ldots, x_{n+1}$ and formulas  (\ref{l1})-(\ref{l3}) hold for $j=0,\ldots,n$.

Consider an interval in $X$ 
$$x=x_n+\frac{\tau-\tau_n}{\tau_{n+1}-\tau_n}(x_{n+1}-x_n),\quad \tau\in[\tau_n,\tau_{n+1}].$$ Then by the induction assumption $\|x_{n+1}-x_n\|_X\le \tau_{n+1}-\tau_n$ it follows that
$$\|x-x_n\|_X\le \tau-\tau_n$$
and
$$\|x-x_0\|_X\le \|x-x_n\|_X+\sum_{k=1}^n\|x_k-x_{k-1}\|_X\le \tau-\tau_0<r.$$
From  hypothesis H2 we obtain
$$\|\Phi'(x)\|\le\ph'(\tau).$$

Consequently we have
\begin{align}\|&\Phi(x_{n+1})-\Phi(x_n)\|_Y\nonumber\\
&=\Big\|\int_0^1\Phi'\big(x_n+t(x_{n+1}-x_n)\big)(x_{n+1}-x_n)dt\Big\|_Y
\le\int_{\tau_n}^{\tau_{n+1}}\ph'(\tau)d\tau\nonumber\\
&=\ph(\tau_{n+1})-\ph(\tau_n)=\psi(\tau_{n+2})-\psi(\tau_{n+1}).\nonumber\end{align}

Observing that $\Phi(x_n)=\Psi(x_{n+1})$ we yield  $$\|\Phi(x_{n+1})-\Psi(x_{n+1})\|_Y\le\psi(\tau_{n+2})-\psi(\tau_{n+1}) .$$

So we are ready to use hypothesis H1 and definition \ref{csf} with
$$\tau'=\tau_{n+1},\quad \tau''=\tau_{n+2},\quad x'=x_{n+1}$$
to obtain $x_{n+2}$ such that $\Psi(x_{n+2})=\Phi(x_{n+1})$ and $$\|x_{n+2}-x_{n+1}\|_X\le \tau_{n+2}-\tau_{n+1}.$$

From the induction assumption we also have 
$$\|x_{n+2}-x_0\|_Y\le \tau_{n+2}-\tau_0.$$
The Lemma is proved. 

Now we are ready to prove the theorem. Since a series $\sum_{k=1}^\infty(\tau_k-\tau_{k-1})$ is convergent then by (\ref{l2}) the series $\sum_{k=1}^\infty(x_k-x_{k-1})$ is also convergent. Therefore the sequence $\{x_k\}$ converges: $x_k\to x_*$. From (\ref{l1}) it follows that $\|x_*-x_0\|_X\le \tau_*-\tau_0<r.$ To finish the proof it remains to pass to the limit in (\ref{l3}).

The Theorem is proved.

\subsection*{Acknowledgement} The author wishes to thank Professor A. V. Arutyunov for useful discussions.

 \end{document}